\documentclass[12pt]{article}
\usepackage{amsbsy,amsfonts,amsmath,amssymb,amsthm,euscript,enumerate}
\usepackage{centernot,graphicx,color,authblk}

\allowdisplaybreaks[1]

\theoremstyle{plain}
\numberwithin{equation}{section}
\newtheorem{thm}{Theorem}[section]

\newtheorem{lmm}[thm]{Lemma}

\newcommand{\conn}{\longrightarrow}

\newcommand{\dc}{d_\mathrm{c}}

\newcommand{\lbeq}[1]{\label{eq:#1}}
\newcommand{\mL}{{\mathbb L}}

\newcommand{\mP}{{\mathbb P}}
\newcommand{\N}{{\mathbb N}}
\newcommand{\nn}{\nonumber}
\newcommand{\opbond}[2]{[#1,#2\rangle}
\newcommand{\pc}{p_\mathrm{c}}
\newcommand{\Proof}[1]{\paragraph{\it #1}}
\newcommand{\QED}{\hspace*{\fill}\rule{7pt}{7pt}\smallskip}
\newcommand{\refeq}[1]{(\ref{eq:#1})}

\newcommand{\vep}{\varepsilon}

\newcommand{\Z}{\mathbb{Z}}
\newcommand{\Zd}{\Z^d}

\title{Hyperscaling for oriented percolation in $1+1$ space-time dimensions}
\author{
Akira Sakai\footnote{{\tt http://www.math.sci.hokudai.ac.jp/\~{}sakai/}}}
\affil{Department of Mathematics, Hokkaido University}
\date{February 27, 2018}

\begin{document}
\maketitle

\begin{abstract}
Consider nearest-neighbor oriented percolation in $d+1$ space-time dimensions.  
Let $\rho,\eta,\nu$ be the critical exponents for the survival probability up 
to time $t$, the expected number of vertices at time $t$ connected from the
space-time origin, and the gyration radius of those vertices, respectively.  We 
prove that the hyperscaling inequality $d\nu\ge\eta+2\rho$, which holds for all 
$d\ge1$ and is a strict inequality above the upper-critical dimension 4, 
becomes an equality for $d=1$, i.e., $\nu=\eta+2\rho$, provided existence of 
at least two among $\rho,\eta,\nu$.  The key to the proof is the recent result 
on the critical box-crossing property by Duminil-Copin, Tassion and 
Teixeira \cite{dctt17}.
\end{abstract}

\section{Introduction and the main results}
Oriented percolation is a time-oriented model of percolation.  It is also 
considered as a discrete-time model for the spread of an infectious disease, 
known as the contact process or the SIS model.  Since it became known to 
exhibit a phase transition and critical behavior, there have been intensive 
researches in both theory and applications in various fields.  Recently, 
a possible association to the laminar-turbulent flow transition was reported 
in \cite{st16}.

Consider the following nearest-neighbor bond oriented percolation on the 
space-time lattice $\mL^d\equiv\{(x,t)\in\Zd\times\Z_+:\|x\|_1+t$ is even\}.  
A pair of vertices $\opbond{(x,s)}{(y,t)}$ is called a bond if $\|x-y\|_1=1$ 
and $t=s+1$.  Each bond $\opbond{(x,t)}{(y,t+1)}$ is either occupied with 
probability $p\in[0,1]$ or vacant with probability $1-p$, independently of the 
other bonds.  Let $\mP_p$ be the associated probability measure.  We say that 
$(x,s)\in\mL^d$ is connected to $(y,t)\in\mL^d$, denoted by $(x,s)\conn(y,t)$, 
if either $(x,s)=(y,t)$ or there is a sequence of occupied bonds 
$\{\opbond{(v_j,j)}{(v_{j+1},j+1)}\}_{j=s}^{t-1}$ from $v_s=x$ to $v_t=y$.  We 
simply write $(x,s)\conn t$ for the event $\bigcup_y\{(x,s)\conn(y,t)\}$, and 
$s\conn(y,t)$ for the event $\bigcup_x\{(x,s)\conn(y,t)\}$.

The major quantities we are interested in are the following.  The first
quantity is the survival probability up to time $t$, defined as
\begin{align}
\theta_t=\mP_p\big((o,0)\conn t\big),
\end{align}
where, and in the rest of the paper, the $p$-dependence is suppressed for
lighter notation.  Since $\{\theta_t\}_{t\in\N}$ is a decreasing sequence of
increasing and continuous functions in $p$, the limit
$\theta_\infty\equiv\lim_{t\uparrow\infty}\theta_t$ is nondecreasing and
right-continuous in $p$.  Let
\begin{align}
\pc=\inf\{p\in[0,1]:\theta_\infty>0\}.
\end{align}
It is proven in \cite{gh02} that $\theta_\infty$ is also left-continuous
in $p$.  In particular, $\theta_\infty=0$ at $p=\pc$, which has not been
proven yet for unoriented percolation in full generality.

The second and third quantities are the expected number of vertices at time
$t$ connected from the origin $(o,0)$ and the gyration radius of those vertices,
defined as
\begin{align}
\chi_t=\sum_x\tau(x,t),&&
\xi_t=\bigg(\frac1{\chi_t}\sum_x|x|^2\tau(x,t)\bigg)^{1/2},
\end{align}
where $\tau(x,t)$ is the two-point function:
\begin{align}
\tau(x,t)=\mP_p\big((o,0)\conn(x,t)\big).
\end{align}
It is first proven in \cite{ab87}, and recently reproved in a much simpler way
in \cite{dct16}, that the critical point is unique in the sense that
\begin{align}
\pc=\sup\bigg\{p\in[0,1]:\sum_{t=0}^\infty\chi_t<\infty\bigg\}.
\end{align}
The sum $\sum_t\chi_t$ is often called the susceptibility.

Now we briefly summarize the basic properties of those quantities readily 
obtained from the definition.  First we note that, by the Markov property 
and translation invariance,
\begin{align}
\theta_{s+t}\ge\theta_s\theta_t,&&
\chi_{s+t}\le\chi_s\chi_t.
\end{align}
With the help of the trivial inequality
$\theta_t\le\chi_t\le(2t+1)^d\theta_t$, we can conclude that there is a 
common relaxation time $\zeta\in[0,\infty]$ such that
\begin{align}\lbeq{relax}
\zeta=\lim_{t\uparrow\infty}\frac{-t}{\log\theta_t}
 =\sup_{t\in\N}\frac{-t}{\log\theta_t}
 =\lim_{t\uparrow\infty}\frac{-t}{\log\chi_t}
 =\inf_{t\in\N}\frac{-t}{\log\chi_t}.
\end{align}
Using the second and forth equalities, we can say that $\zeta$ is bounded away 
from zero and infinity when $p<\pc$, implying exponential decay of $\theta_t$ 
and $\chi_t$ in $t$ in the subcritical regime.  This is not the case at the 
critical point.  Moreover, $\chi_t$ is nondecreasing in $t$ at $p=\pc$, because, 
otherwise, there must be a $t_0\in\N$ such that $\chi_{t_0}<1$, which together 
with submultiplicativity implies exponential decay of $\chi_t$ and convergence 
of the susceptibility $\sum_t\chi_t$ at $p=\pc$, such as
\begin{align}
\sum_{t=0}^\infty\chi_t=\sum_{n=0}^\infty\sum_{k=0}^{t_0-1}\chi_{nt_0+k}\le
 \sum_{n=0}^\infty\chi_{t_0}^n\sum_{k=0}^{t_0-1}\chi_k<\infty,
\end{align}
which is a contradiction to the result in \cite{an84}: $\sum_t\chi_t=\infty$ at 
$p=\pc$.

Let $\rho,\eta,\nu$ be the critical exponents for the above quantities at 
$p=\pc$: as $t\uparrow\infty$,
\begin{align}
\theta_t\approx t^{-\rho},&&
\chi_t\approx t^\eta,&&
\xi_t\approx t^\nu,
\end{align}
where $f\approx g$ means that $(\log f)/\log g$ goes to 1 in the prescribed 
limit, allowing corrections of slowly varying functions.  In higher dimensions 
$d\gg4$ ($d>4$ is enough for sufficiently spread-out models), the lace 
expansion converges and the above critical exponents take on their mean-field 
values $\rho=1$, $\eta=0$ and $\nu=1/2$: the values for branching random walk 
\cite{cs08,cs11,hh13,hs02,ny93,ny95,s01}.  In lower dimensions, on the other 
hand, only numerical values and predictions due to non-rigorous 
renormalization-group methods are available (see Table~\ref{tab:crexp}).

\begin{table}[b]
\begin{center}
\caption{\label{tab:crexp}Predicted values of the critical exponents in various 
dimensions (e.g., \cite{o04}).}
\medskip
{\renewcommand\arraystretch{1.3}
\begin{tabular}{|c||c|c|c|c|c|}
\hline
 & $d=1$ & $d=2$ & $d=3$ & $d=4-\vep$ & $d\ge4$ \\
\hline
$\rho$ & 0.159464 & 0.451 & 0.73 & $1-\frac14\vep-0.01283\vep^2$ & 1 \\
\hline
$\eta$ & 0.313686 & 0.230 & 0.12 & $\frac1{12}\vep+0.03751\vep^2$ & 0 \\
\hline
$\nu$ & 0.632613 & 0.568 & 0.526 & $\frac12+\frac1{48}\vep+0.008171\vep^2$ & $\frac12$ \\
\hline
$\gamma$ & 2.277730 & 1.60 & 1.25 & $1+\frac16\vep+0.06683\vep^2$ & 1 \\
\hline
$\mu$ & 1.733847 & 1.295 & 1.105 & $1+\frac1{12}\vep+0.02238\vep^2$ & 1 \\
\hline
\end{tabular}}
\end{center}
\end{table}

In this paper, we prove the following theorem.

\begin{thm}\label{thm:main}
\begin{enumerate}[(i)]
\item
For any $d\ge1$, $p\in[0,1]$ and $t\in\N$, we have
\begin{align}\lbeq{hypineq1}
\chi_t\le\frac43(4\xi_t+1)^d\,\theta_{t/2}^2,
\end{align}
which implies the hyperscaling inequality (assuming existence of 
$\rho,\eta,\nu$)
\begin{align}\lbeq{hypineq2}
d\nu\ge\eta+2\rho.
\end{align}
\item
Let $d=1$ and $p=\pc$.  Then, there is a $K>0$ such that, for any $t\in\N$,
\begin{align}\lbeq{hypeq1}
\chi_t\ge K\xi_t\theta_t^2,
\end{align}
which implies the hyperscaling equality (assuming existence of 
at least two among $\rho,\eta,\nu$)
\begin{align}\lbeq{hypeq2}
\nu=\eta+2\rho.
\end{align}
\end{enumerate}
\end{thm}

\paragraph{Remark:}
\begin{enumerate}
\item
The inequality \refeq{hypineq1} was first derived in \cite{s02}.  Since its 
proof is easy and short, we will show it again for convenience.  It was used 
in \cite{s02} to prove two other hyperscaling inequalities that also involve 
critical exponents defined in the off-critical regime.  For example, if the 
susceptibility $\sum_t\chi_t$ and the relaxation time $\zeta$ diverge as 
$p\uparrow\pc$ as $(\pc-p)^{-\gamma}$ and $(\pc-p)^{-\mu}$ respectively, 
then, for any $d\ge1$, we have
\begin{align}\lbeq{hypineq3}
(d\nu-2\rho+1)\mu\ge\gamma.
\end{align}
If we replace those critical exponents in \refeq{hypineq2} and \refeq{hypineq3} 
by their mean-field values, then we obtain $d\ge4$, which is a complement to 
the aforementioned lace-expansion results.  Therefore, the upper-critical 
dimension $\dc$ for oriented percolation is 4.
\item
In general, hyperscaling inequalities are believed to be equalities below and 
at the model-dependent upper-critical dimension.  The values in 
Table~\ref{tab:crexp} seem to support this belief.  The identity \refeq{hypeq2} 
proves that it is indeed the case for at least $d=1$.  For unoriented 
percolation, for which $\dc=6$, similar results are proven in 2 dimensions by 
Kesten~\cite{k87} using the Russo-Seymour-Welsh theorem on the critical 
box-crossing property \cite{r78,r81,sw78}.  Since the known critical exponents 
for 2-dimensional unoriented percolation are rational numbers (e.g., 
$\beta=5/36$ and $\gamma=43/18$), it is natural to believe that there must 
be some balance (i.e., hyperscaling equalities) among those critical 
exponents.  On the other hand, since the values in Table~\ref{tab:crexp} do 
not seem to be rational numbers, the hyperscaling equality \refeq{hypeq2} is 
even more surprising.
\item
The main reason why the right-hand side of \refeq{hypineq1} is bigger than its 
left-hand side is due to the inequality
\begin{align}\lbeq{keyineq}
\tau(x,t)&=\mP_p\big((o,0)\conn(x,t)\big)\nn\\
&\le\mP_p\big((o,0)\conn t/2,~t/2\conn(x,t)\big)=\theta_{t/2}^2,
\end{align}
where, and in the rest of the paper, we do not care much about possibilities 
of, e.g., $t/2$ not being an integer, since it is easy (but cumbersome) to make 
the argument rigorous if we introduce floor functions, etc.  The last equality 
in \refeq{keyineq} is due to reversibility: if we change the direction 
of each bond and redefine the connectivity in the time-decreasing direction, 
then we have the identity $\mP_p(t/2\conn(x,t))=\theta_{t/2}$.  
\item
The following theorem on the critical box-crossing property is the key to show 
the opposite inequality to \refeq{keyineq}:

\begin{thm}[Theorem~1.3 in \cite{dctt17}]\label{thm:dctt17}
Let
\begin{align}
V_p(w,t)&=\mP_p\Big([0,w]\times[0,t]\text{ is crossed vertically}\Big),\\
H_p(w,t)&=\mP_p\Big([0,w]\times[0,t]\text{ is crossed from left to right}\Big).
\end{align}
There exist a constant $\vep\in(0,1)$ and an increasing sequence of integers 
$\{w_t\}_{t\in\N}$ such that, for all $t\in\N$,
\begin{gather}
\vep\le V_{\pc}(w_t,3t)\le V_{\pc}(3w_t,t)\le1-\vep,\lbeq{vcrossing}\\[5pt]
\vep\le H_{\pc}(3w_t,t)\le H_{\pc}(w_t,3t)\le1-\vep.\lbeq{hcrossing}
\end{gather}
\end{thm}

\bigskip

We will also use \refeq{vcrossing}--\refeq{hcrossing} to control an upper 
bound on $\tau(x,t)$ for $x>jw_t$ that decays exponentially in $j\in\N$ 
(see Lemma~\ref{lmm:main} below).  
This is a key element to show that $w_t$ is bounded below by an 
$\vep$-dependent positive multiple of $\xi_t$.  
\item
Applying \refeq{hypineq1} and \refeq{hypeq1} to \cite[(5.1)]{s02} and its 
reverse, respectively, we can readily show that the 
hyperscaling inequality \refeq{hypineq3} also becomes an equality for $d=1$, 
i.e.,
\begin{align}
(\nu-2\rho+1)\mu=\gamma.
\end{align}
\item
It is easy to show that the hyperscaling inequality \refeq{hypineq2} holds for 
other finite-range models of oriented percolation and the contact process.  
It should not be so difficult to prove Theorem~\ref{thm:dctt17} for the 
nearest-neighbor models of oriented site percolation and the contact process, 
hence the hyperscaling equality \refeq{hypeq2} for $d=1$.  However, it is not 
so obvious to prove a similar statement to Theorem~\ref{thm:dctt17} for 
longer-range models.  This may be worth further investigation.
\end{enumerate}

\section{Proof of Theorem~\ref{thm:main}}
\Proof{Proof of Theorem~\ref{thm:main}(i).}
It suffices to prove the inequality \refeq{hypineq1}, as the hyperscaling 
inequality \refeq{hypineq2} immediately follows by using \refeq{hypineq1} 
at $p=\pc$ (and assuming existence of the three critical exponents).  First 
we note that
\begin{align}
\chi_t=\frac1{\xi_t^2}\sum_x|x|^2\tau(x,t)\ge4\sum_{x:|x|\ge2\xi_t}\tau(x,t),
\end{align}
hence
\begin{align}
\frac34\chi_t\le\sum_{x:|x|\le2\xi_t}\tau(x,t).
\end{align}
By \refeq{keyineq}, the right-hand side is further bounded by 
$(4\xi_t+1)^d\theta_{t/2}^2$.  This completes the proof of \refeq{hypineq1}.
\QED

\bigskip

To prove Theorem~\ref{thm:main}(ii), we first assume the following key lemma:

\begin{lmm}\label{lmm:main}
Let $d=1$ and $p=\pc$.  Let $\vep\in(0,1)$ and $w_t$ be the same as in 
Theorem~\ref{thm:dctt17}. 
\begin{enumerate}[(i)]
\item
For any $t\in\N$ and any $x\in[-\frac12w_t,\frac12w_t]$,
\begin{align}\lbeq{lmm1}
\tau(x,t)\ge\vep^6\theta_t^2.
\end{align}
\item
For any $j,t,x\in\N$ with $j\ge2$ and $jw_t<x\le(j+1)w_t$, 
\begin{align}\lbeq{lmm2}
\tau(x,t)\le\vep^{-4}\theta_t^2(1-\vep)^{j-2}.
\end{align}
\end{enumerate}
\end{lmm}

\Proof{Proof of Theorem~\ref{thm:main}(ii) assuming Lemma~\ref{lmm:main}.}
Again, it suffices to prove the inequality \refeq{hypeq1}, as the equality 
\refeq{hypeq2} is a result of the hyperscaling inequality \refeq{hypineq2} for 
$d=1$ and the opposite inequality $\nu\le\eta+2\rho$ that immediately follows 
from \refeq{hypeq1}.

To prove \refeq{hypeq1}, we first note that, by \refeq{lmm1}, 
\begin{align}\lbeq{chitlbd1}
\chi_t\ge2\sum_{x=1}^{\frac12w_t}\tau(x,t)\ge\vep^6w_t\theta_t^2.
\end{align}
To complete the proof, it suffices to show that $w_t$ is bounded below by 
a positive multiple of $\xi_t$.  However, by definition,
\begin{align}
\xi_t^2=2\sum_{x=1}^\infty x^2\frac{\tau(x,t)}{\chi_t}&=2\bigg(\sum_{x=1}^{2w_t}
 x^2\frac{\tau(x,t)}{\chi_t}+\sum_{j=2}^\infty\sum_{x=jw_t+1}^{(j+1)w_t}x^2
 \frac{\tau(x,t)}{\chi_t}\bigg)\nn\\
&\le2w_t^2\bigg(4+\sum_{j=2}^\infty
 (j+1)^2\sum_{x=jw_t+1}^{(j+1)w_t}\frac{\tau(x,t)}{\chi_t}\bigg).
\end{align}
Then, by using \refeq{lmm2}--\refeq{chitlbd1}, we obtain
\begin{align}
\xi_t^2&\stackrel{\refeq{chitlbd1}}\le2w_t^2\bigg(4+\frac1{\vep^6\theta_t^2}
 \sum_{j=2}^\infty(j+1)^2\max_{jw_t<x\le(j+1)w_t}\tau(x,t)\bigg)\nn\\
&\stackrel{\refeq{lmm2}}\le2w_t^2\bigg(4+\vep^{-10}\sum_{j=2}^\infty(j+1)^2
 (1-\vep)^{j-2}\bigg).
\end{align}
As a result,
\begin{align}
\chi_t\ge\underbrace{\frac{\vep^6}{\sqrt2}\bigg(4+\vep^{-10}\sum_{j=2}^\infty
 (j+1)^2(1-\vep)^{j-2}\bigg)^{-1/2}}_{=K}\xi_t\theta_t^2.
\end{align}
This completes the proof of \refeq{hypeq1}.
\QED

\bigskip

The rest of the paper is devoted to showing Lemma~\ref{lmm:main}.

\Proof{Proof of Lemma~\ref{lmm:main}(i).}
First we note that, for $1\le x\le\frac12w_t$, the event $(o,0)\conn(x,t)$ 
occurs if the following four increasing events occur:
\begin{itemize}
\item
$(o,0)\conn t$ in $[-w_t,w_t]\times[0,t]$,
\item
$0\conn(x,t)$ in $[x-w_t,x+w_t]\times[0,t]$,
\item
$[-\frac32w_t,\frac32w_t]\times[0,t]$ is crossed from left to right,
\item
$[-\frac32w_t,\frac32w_t]\times[0,t]$ is crossed from right to left.
\end{itemize}
The last two events take care of the possibility that the forward cluster from 
the origin $(o,0)$ and the backward cluster from $(x,t)$ do not collide.  Using 
the FKG inequality (see, e.g., \cite{g99}), translation invariance and the 
reversibility explained below \refeq{keyineq}, we obtain
\begin{align}
\tau(x,t)\ge\mP_p\Big((o,0)\conn t\text{ in }[-w_t,w_t]\times[0,t]\Big)^2
  H_p(3w_t,t)^2.
\end{align}
We further note that the event $(o,0)\conn t$ in $[-w_t,w_t]\times[0,t]$ occurs 
if the following three increasing events occur:
\begin{itemize}
\item
$(o,0)\conn t$, 
\item
$[0,w_t]\times[0,t]$ is crossed vertically,
\item
$[-w_t,0]\times[0,t]$ is crossed vertically.
\end{itemize}
Again, by the FKG inequality, translation invariance and the 
monotonicity $V_p(w_t,t)\ge V_p(w_t,3t)$, we obtain
\begin{align}
\mP_p\Big((o,0)\conn t\text{ in }[-w_t,w_t]\times[0,t]\Big)\ge\theta_tV_p(w_t,
 3t)^2,
\end{align}
hence
\begin{align}
\tau(x,t)\ge\theta_t^2V_p(w_t,3t)^4H_p(3w_t,t)^2.
\end{align}
The inequality \refeq{lmm1} follows from the above inequality at 
$p=\pc$ and \refeq{vcrossing}--\refeq{hcrossing}.
\QED

\Proof{Proof of Lemma~\ref{lmm:main}(ii).}
Recall that $j\ge2$ and $jw_t<x\le(j+1)w_t$.  If $(o,0)\conn(x,t)$, then the 
following three independent events occur:
\begin{itemize}
\item
$(o,0)$ is connected to the boundary $\partial B_o$ of the box 
$B_o\equiv[-w_t,w_t]\times[0,t]$,
\item
$[w_t,(j-1)w_t]\times[0,t]$ is crossed from left to right, 
\item
$(x,t)$ is connected from the boundary $\partial B_x$ of the box 
$B_x=[(j-1)w_t,(j+2)w_t]\times[0,t]$.
\end{itemize}
By this observation and using 
$H_p((j-2)w_t,t)\le H_p(w_t,t)^{j-2}\le H_p(w_t,3t)^{j-2}$, we obtain 
\begin{align}
\tau(x,t)\le\mP_p\Big((o,0)\conn\partial B_o\Big)H_p(w_t,3t)^{j-2}
 \mP_p\Big(\partial B_x\conn(x,t)\Big).
\end{align}
However, by reversibility and monotonicity, we have 
\begin{align}
\mP_p\Big(\partial B_x\conn(x,t)\Big)\le\mP_p\Big((o,0)\conn\partial B_o\Big).
\end{align}
Therefore,
\begin{align}
\tau(x,t)\le\mP_p\Big((o,0)\conn\partial B_o\Big)^2H_p(w_t,t)^{j-2}.
\end{align}

To bound the probability on the right-hand side by $\theta_t$, we borrow the 
idea in the proof of \cite[(4.7)]{dctt17}.  First, we note that $(o,0)\conn t$ 
if the following three increasing events occur:
\begin{itemize}
\item
$(o,0)\conn\partial B_o$,
\item
$[0,w_t]\times[0,t]$ is crossed vertically, 
\item
$[-w_t,0]\times[0,t]$ is crossed vertically.
\end{itemize}
By the FKG inequality, translation invariance and 
the monotonicity $V_p(w_t,t)\ge V_p(w_t,3t)$, we obtain
\begin{align}\lbeq{keyidea2}
\theta_t\ge\mP_p\Big((o,0)\conn\partial B_o\Big)V_p(w_t,3t)^2.
\end{align}

To summarize the above computations at $p=\pc$, we arrived at
\begin{align}
\tau(x,t)\le\bigg(\frac{\theta_t}{V_{\pc}(w_t,3t)^2}\bigg)^2H_{\pc}(w_t,
 3t)^{j-2}\le\vep^{-4}\theta_t^2(1-\vep)^{j-2},
\end{align}
as required.
\QED

\section*{Acknowledgements}
This work was initiated when I started preparation for the Summer School in 
Mathematical Physics, held at the University of Tokyo from August 25 through 
27, 2017.  I am grateful to the organizers, Yasuyuki Kawahigashi and Yoshiko 
Ogata, for the opportunity to speak at the summer school and meet with many 
researchers in the laminar-turbulent flow transition.  Finally, I would like to 
thank Alessandro Giuliani for his support during the refereeing process and 
an anonymous referee for valuable comments to the earlier version to this 
paper.


\begin{thebibliography}{99}
\bibitem{ab87}M. Aizenman and D.J. Barsky.
\newblock Sharpness of the phase transition in percolation models.
\newblock \emph{Commun. Math. Phys.}~\textbf{108} (1987): 489--526.

\bibitem{an84}M. Aizenman and C.M. Newman.
\newblock Tree graph inequalities and critical behavior in percolation models.
\newblock \emph{J. Stat. Phys.}~\textbf{36} (1984): 107--143.

\bibitem{cs08}L.-C. Chen and A. Sakai.
\newblock Critical behavior and the limit distribution for long-range oriented
percolation.~I.
\newblock \emph{Probab. Theory Related Fields}~\textbf{142} (2008): 151--188.

\bibitem{cs11}L.-C. Chen and A. Sakai.
\newblock Asymptotic behavior of the gyration radius for long-range
self-avoiding walk and long-range oriented percolation.
\newblock \emph{Ann. Prob.}~\textbf{39} (2011): 507--548.

\bibitem{dct16}H. Duminil-Copin and V. Tassion.
\newblock A new proof of the sharpness of the phase transition for Bernoulli
percolation and the Ising model.
\newblock \emph{Commun. Math. Phys.}~\textbf{343} (2016): 725--745.

\bibitem{dctt17}H. Duminil-Copin, V. Tassion and A. Teixeira.
\newblock The box-crossing property for critical two-dimensional oriented
percolation.
\newblock To appear in \emph{Probab. Theory Related Fields}.
\texttt{arXiv:1610.10018}.

\bibitem{g99}G. Grimmett.
\newblock \emph{Percolation} (2nd ed., Springer, 1999).

\bibitem{gh02}G. Grimmett and P. Hiemer.
\newblock Directed percolation and random walk.
\newblock \emph{In and Out of Equilibrium} (V. Sidoravicius ed.,
Birkh\"auser, 2002): 273--297.

\bibitem{hh13}R. van der Hofstad and M. Holmes.
\newblock The survival probability and $r$-point functions in high dimensions.
\newblock \emph{Ann. Math.}~\textbf{178} (2013): 665--685.

\bibitem{hs02}R. van der Hofstad and G. Slade.
\newblock A generalised inductive approach to the lace expansion.
\newblock \emph{Probab. Theory Related Fields}~\textbf{122} (2002): 389--430.

\bibitem{k87}H. Kesten.
\newblock Scaling relations for $2D$-percolation.
\newblock \emph{Commun. Math. Phys.}~\textbf{109} (1987): 109--156.

\bibitem{ny93}B.G. Nguyen and W.-S. Yang.
\newblock Triangle condition for oriented percolation in high dimensions.
\newblock \emph{Ann. Prob.}~\textbf{21} (1993): 1809--1844.

\bibitem{ny95}B.G. Nguyen and W.-S. Yang.
\newblock Gaussian limit for critical oriented percolation in high dimensions.
\newblock \emph{J. Stat. Phys.}~\textbf{78} (1995): 841--876.

\bibitem{o04}G. \'Odor.
\newblock Universality classes in nonequilibrium lattice systems.
\newblock \emph{Rev. Mod. Phys.}~\textbf{76} (2004): 663--724.

\bibitem{r78}L. Russo
\newblock A note on percolation.
\newblock \emph{Z. Wahrscheinlichkeitstheor. verw. Geb.}~\textbf{43} (1978): 
39--48.

\bibitem{r81}L. Russo.
\newblock On the critical percolation probabilities.
\newblock \emph{Z. Wahrscheinlichkeitstheor. verw. Geb.}~\textbf{56} (1981): 
229--237.

\bibitem{s01}A. Sakai.
\newblock Mean-field critical behavior for the contact process.
\newblock \emph{J. Stat. Phys.}~\textbf{104} (2001): 111--143.

\bibitem{s02}A. Sakai.
\newblock Hyperscaling inequalities for the contact process and oriented
percolation.
\newblock \emph{J. Stat. Phys.}~\textbf{106} (2002): 201--211.

\bibitem{st16}M. Sano and K. Tamai.
\newblock A universal transition to turbulence in channel flow.
\newblock \emph{Nature Phys.}~\textbf{12} (2016): 249--253.

\bibitem{sw78}P.D. Seymour and D.J.A. Welsh.
\newblock Percolation probabilities on the square lattice. 
\newblock \emph{Ann. Discrete Math.}~\textbf{3} (1978): 227--245.

\end{thebibliography}
\end{document}